\documentclass[12pt,reqno]{amsart}
\usepackage{amsmath}%
\usepackage{amsfonts}%
\usepackage{amssymb}%
\usepackage[margin=.5in]{geometry}
\usepackage{subfig}
\usepackage{array,multirow}
\usepackage{ifpdf}
\ifpdf
\usepackage{graphicx}
\else
\usepackage[draft]{graphicx}
\fi
\usepackage{xcolor}
\usepackage{enumerate}
\usepackage[T1]{fontenc}
\usepackage{caption}
\usepackage{booktabs,multirow}
\usepackage{enumerate}
\usepackage{listings}
\usepackage{caption}
\usepackage{hyperref}
\usepackage{caption}
\usepackage{threeparttable}
\usepackage{bm}
\usepackage{csquotes}
\usepackage{float}

\MakeOuterQuote{"}

\hypersetup{
    colorlinks=true,
    linkcolor=blue,
    filecolor=magenta,      
    urlcolor=cyan,
    pdfpagemode=UseNone,pdfstartview=
    }

\definecolor{dkgreen}{rgb}{0,0.6,0}
\definecolor{gray}{rgb}{0.5,0.5,0.5}
\definecolor{mauve}{rgb}{0.58,0,0.82}

\AtBeginDocument{}

\lstdefinestyle{myMathematica}{
frame=tb,
  language=Mathematica,
  aboveskip=3mm,
  belowskip=3mm,
  showstringspaces=false,
  columns=flexible,
  basicstyle={\ttfamily},
  numbers=none,
  numberstyle=\tiny\color{gray},
  keywordstyle=\color{blue},
  commentstyle=\color{mauve},
  stringstyle=\color{dkgreen},
  breaklines=true,
  breakatwhitespace=true,
	morekeywords={AsymptoticSolve},
  tabsize=3
 }
\begin{document}

\title[]{Determining radii of convergence of fractional power expansions around singular points of algebraic functions}

\author{Dominic C. Milioto}
\email{icorone@hotmail.com}

\date{\today}
\subjclass[2010]{Primary 1401; Secondary 1404} %
\begin{abstract}
The purpose of this paper is to introduce the branching geometry of algebraic functions around singular points and to describe a simple method of determining radii of convergence of their power expansions in terms of those singular points.  Branching geometries are categorized into six types.  Then a method is presented to determine radii of convergence of branch power expansions using analytic continuation and the identification of convergence-limiting singular points.  Test cases exhibiting a variety of branching morphologies are analyzed, and convergence results obtained through analytic continuation are checked against Root Tests of the associated power series. All Root Tests agreed well with the results obtained by analytic continuation.  Mathematica ver. 12.3 was used to implement the numeric algorithms.  

\end{abstract}
\maketitle
\section{Introduction}

This paper is about the the branching geometry and power expansions of an algebraic function $w(z)$ defined implicitly by the irreducible $n$-degree expression in $w$:
\begin{equation}
f(z,w)=a_0(z)+a_1(z)w+a_2(z)w^2+\cdots+a_{n}(z)w^{n}=0
\label{eqn:eqn001}
\end{equation}
with $z$ and $w$ complex variables and the coefficients, $a_i(z)$, polynomials in $z$ with rational coefficients.  By the Implicit Function Theorem, (\ref{eqn:eqn001}) defines locally, an analytic function $w(z)$ when $\displaystyle \frac{\partial f}{\partial w}\neq 0$.  The solution set of (\ref{eqn:eqn001}) defines an algebraic curve, $w(z)$, and it is known from the general theory of algebraic functions that $w(z)$ can be described in a disk centered at $z_0$ by $n$ fractional power series called Puiseux series with a radii of convergence extending at least to the distance to the nearest singular point.  This paper develops a method to determine the exact radius of convergence of each series in terms of the singular points of $w(z)$.  Puiseux series can be computed by the method of Newton Polygon\cite{Milioto}.  Both custom code designed by this author and Mathematica's $\texttt{AsymptoticSolve}$ function were used to compute Puiseux expansions for this study.
 
 $w(z)$ for $n>1$ is multi-valued with a helical branching geometry not unlike a composition of $k$-root functions albeit more complex.  A Puiseux expansion of $w(z)$ at a point $z_0$ is a set of $n$ fractional power expansion in $z$ given by 
\begin{equation}
 P_i(z)=w_i(z)=\sum_{k=r}^\infty a_k (z-z_0)^{\frac{m_k}{c}},\quad i=1,2,\cdots,n
\label{eqn:eqn002}
\end{equation} 
where $z$ lies in the domain of convergence of the series and the terms in the expression defined below.  In this paper, Puiseux series are written as $\displaystyle P_i(z)=\sum_{k=r}^\infty a_k z^{\frac{m_k}{c}}$  with $z^{\frac{m_k}{c}}$  interpreted as principal-value and a specific $t$ value of the series given by $P_i(t-z_0)$.  See Appendix \ref{appendix:appB} for more information about Puiseux series and the format used here.     

  The derivative of $w(z)$ at a point $(z,w)=(p,q)$ can be computed as
\begin{equation}
\frac{dw}{dz}=\lim_{(z,w)\to(p,q)} \left(-\frac{f_z(z,w)}{f_w(z,w)}\right)
\label{eqn:eqn003}
\end{equation}

 when this limit exist. A singularity $\{z_s,w_s\}$  of $w(z)$ is a point where the limit does not exist and in this paper, the term \enquote{singular point} refers to the $z$-component of the singularity.

 Singular points can be computed numerically by solving for the zeros of the resolvent $R(f,f_w)$. 

  In the case of an indeterminate form $\displaystyle \frac{0}{0}$ of (\ref{eqn:eqn003}), $w(z)$ may have a removable singularity at $z_s$ if the limit exist.  This can be determined by computing the derivatives of the associated Puiseux expansions and identifying $1$-cycle series having finite derivatives at $0$ (recall substitution of $(t-z_s)$):  If $f(z_s,w_0)=f_z(z_s,w_0)=f_w(z_s,w_0)=0$ , then for each $1$-cycle series for which $\displaystyle\frac{d P_i}{dz}(0)$ is finite,  
\begin{equation}
\lim_{(z,w)\to(z_s,w_0)}\left(-\frac{f_z(z,w_0)}{f_w(z,w_0)}\right)=\frac{dP_i}{dz}\left(0\right);\quad w_0=P_i(0).
\label{eqn:eqn004}
\end{equation} 

Test Case $3$ below includes example calculations with removable singular points.
\section{Conventions used in this paper}
\begin{enumerate}
\item This work studies polynomials given by (\ref{eqn:eqn001}) with random rational coefficients.

\vspace{10pt}\item Although a singularity of $w(z)$ is a point $\{z_s,w_s\}\in\mathbb{C}^2$, \enquote{singular point} in this paper refers to the $z_s$ term.

\vspace{10pt}\item The concept of \enquote{branch} is used throughout this paper and refers to a part of $w(z)$, often multi-valued, and  represented by a set of Puiseux series with finite radii of convergence.  

\vspace{10pt}\item Reference is made to a \enquote{base} singular point $s_b$.  This refers to a center of expansion of Puiseux series with $s_b$ a singular point.  In the procedure described below, the branch surfaces about $s_b$ are  analytically continued over other singular points $s_i$ in order of increasing distance from the base singular point until convergence-limiting singular points (CLSPs) are encountered for each branch sheet.

\vspace{10pt} \item The Puiseux expansions, $\{P_i\}$, of $w(z)$ at a point $s_b$, consists of a set of $n$ fractional power expansions in terms of $\displaystyle z^{1/c}$ where $c$ is a positive integer and can be different for different series in the set.  $c$ is both the cycle size of the series and cycle size of the branch represented by the series.  For example, a power series $\displaystyle \sum_{k=r}^{\infty} a_k z^{\frac{m_k}{3}}$ has a cycle size of $3$ and represents a $3$-cycle branch:  This branch has three coverings over a deleted neighborhood of the expansion center, and the geometry of the branch is similar to the geometry of $\displaystyle z^{1/3}$.  See Figure \ref{figure:figure3001} for an example branch plot. 

\vspace{10pt} \item The Puiseux expansions of $w(z)$ at a point are grouped into conjugate classes. For example, a $5$-cycle branch of $w(z)$ is expanded into five Puiseux series in powers of $z^{1/5}$, one series for each single-valued sheet of the branch.  These five series make up a single $5$-cycle conjugate class.  The sum of the conjugate classes at a point $z$ is always equal to the degree of the function in $w$.  For example, a power expansion of a $10$-degree function consist of the set $\{P_i\}$ such that the sum of the conjugate  types is $10$.  This could consists of a single $10$-cycle conjugate class containing ten series, or three different $3$-cycle conjugate classes and a single $1$-cycle conjugate class or some other combination of conjugate classes adding up to $10$. Each series member in a conjugate class can be generated by conjugation of a member of the class as follows:  

 Let 
\begin{equation}
P_k(z)=\sum_{i=r}^{\infty} a_i z^{\frac{m_i}{c}}
\label{eqn:eqn300}
\end{equation}

 be the $k$-th member of a $c$-cycle conjugate class of Puiseux series where all $\frac{m_i}{c}$ exponents are placed under a least common denominator $c$ and $\displaystyle z^{\frac{m_i}{c}}$ is the principal-valued root as described above.  Then the $c$ members of this conjugate class can be generated via conjugation of (\ref{eqn:eqn300}) as follows:
\begin{equation}
P_j(z)=\sum_{k=r}^{\infty} a_k \left(e^{\frac{2 j\pi i}{c}}\right)^{m_k} z^{\frac{m_k}{c}};\quad j=0,1,\cdots,c-1.
\label{eqn:eqn301}
\end{equation}

\vspace{10pt}\item  If an $n$-degree function expands into a single conjugate class at $s_b$, the underlying branch as well as the power series are said to \enquote{full-ramify} at $s_b$.  For example, if an expansion of $w(z)$ of degree five produces $5$ series in terms of $z^{1/5}$, the function fully-ramifies into a $5$-cycle branch with the five series belonging to the same $5$-cycle conjugate class.  The branch is morphologically similar to $\displaystyle z^{1/5}$.  An expansion of an $n$-degree function can include multiple $k$-cycle conjugate classes.  For example, a $10$-degree function can expand into five sets of $2$-cycle conjugate classes with each class containing two sets of $2$-cycle series.

 \vspace{10pt}\item The Puiseux series generated in this paper are ordered according to the value of each series at a singular point reference $p_r=s_b+1/3 r$ where $r$ is the distance to the nearest singular point and $s_b$ is the expansion center.  The ordering is first by real part then imaginary part of the set $\{P_i(p_r-s_b)\}$.
\
\vspace{10pt}\item This discussion makes use of the phrase, \enquote{extending a branch over a singular point.}  This is in reference to the discussion below about convergence-limiting singular points (CLSPs) and the branches they affect.
\vspace{10pt}\item $R$ is a positive real number representing the radius of convergence of a power series centered at a singular point $s_b$.  The value of $R$ is expressed in terms of the associated CLSP.  For example, if a power expansion has a center at the tenth singular point, $s_{10}$, and its CLSP was found to be $s_{25}$, then $R=|s_{10}-s_{25}|$.  This notation is presented as the exact symbolic expression for radius of convergence.
\vspace{10pt}\item Power series for algebraic branches are fractional power series called Puiseux series.  Puiseux series generated for this work were computed both by Mathematica's $\texttt{AsymptoticSolve}$ and custom code developed by this author.  Note:  $\texttt{AsymptoticSolve}$ generates Puiseux series in powers of $(z-z_0)$ where $z_0$ is the expansion center.

\end{enumerate}

\section{Branch and Series Types}

 The term \enquote{branch} refers to an analytic portion of $w(z)$ in a deleted  disc $D(s_b,R)$ of the complex $z$-plane and is represented by a set of convergent Puiseux series $\{P_i\}$ with center $s_b$.  It is helpful to envision branches as similar to part of the real or imaginary surfaces of a fractional power  $z^{k/c}$ over the complex $z$-plane in the disk $D(0,R)$ albeit often having a more convoluted shape.  Appendix \ref{appendix:appA} are simple examples of branch types covered in this paper. 

 Algebraic branches can be categorized by their cycle type and order of their associated conjugate class of Puiseux series:  The terms of a $c$-cycle series are in the form of $\displaystyle a_k z^{\frac{m_k}{c}}$ with $c$ the least common denominator (LCD) of exponents.  And the order is the smallest non-zero $m_k$ of these exponents.  For example, the exponents of 
$$
z^{3/2}+z^2+z^3+z^{21/3}
$$
in LCD form are $\{9/6,12/6,18/6,42/6\}$.  This is a $6$-cycle of order $9$.  Likewise, exponents of  
$$
 z^{-4}+z^{-5/2}+z^{-3/2}+z+z^{1/3}
$$
in LCD form are $\{-24/6,-15/6,-9/6,6/6,2/6\}$.  This is a $6$-cycle of order $-24$.

 The cycle size and order of a Puiseux expansion determines the geometric shape or \enquote{morphology} of the branch.  Branches are categorized as $X_p^q$ with $p$ the cycle type, $q$ the order, and a letter $X$ identifying the branch morphology according to the classification below.  If the order or cycle size is one, the corresponding subscript or superscript is omitted in the branch descriptor.

\subsection{Branch Types}
\label{sec:bt001}

In the following branch descriptors, all exponents $\displaystyle \frac{q}{p}$ of a series are presumed placed under a least common denominator $p$.

\begin{description}
\item[\textbf{Type \boldsymbol{$T$}}]  Power series with positive integer powers (Taylor series).  These are $1$-cycle branches.
\vspace{5pt}
 \item[\textbf{Type \boldsymbol{$E$}}] $1$-cycle $T$ branch with a removable singular point at its center.
\vspace{5pt}
\item[\textbf{Type \boldsymbol{$F_p^q$}}] $p$-cycle branch with $p>1$ of order $q$ with non-negative exponents and lowest non-zero exponent $\displaystyle \frac{q}{p}$ with $q>p$.  These branches are multi-valued consisting of $p$ single-valued sheets with a finite tangent at the singular point.  An example $F_2^3$ series is $z^{3/2}+z^2+\cdots$.
\vspace{5pt}
\item[\textbf{Type \boldsymbol{$V_p^q$}}] $p$-cycle branch with $p>1$ of order $q$ with non-negative exponents and lowest non-zero exponent $\displaystyle \frac{q}{p}$ with $q<p$ and vertical tangent at center of expansion.  An example $V_4^3$ series is $z^{3/4}+z^2+\cdots$.
\vspace{5pt}
\item[\textbf{Type \boldsymbol{$P_p^q$}}]  $p$-cycle branch unbounded at center with $p>1$ of order $q$ having negative exponents with lowest negative exponent $\displaystyle\frac{q}{p}$.   An example $3$-cycle $P$ series of order $-1$ is $z^{-1/3}+z^2+\cdots$.  An example $3$-cycle $P$ series of order $-5$ is $z^{-5/3}+z^{-1/3}+\cdots$. 
 \vspace{5pt}
\item[\textbf{Type \boldsymbol{$L^q$}}]  Branch with Laurent series of order $q$ as the Puiseux series.  An example $L^{-2}$ series is $1/z^2+1/z+z^2+\cdots$
\end{description}

\section{Brief overview}

 A fractional power expansion of $w(z)$ has a finite radius of convergence and usually represents a small portion of $w(z)$ due to the presence of convergence-limiting singular points (CLSPs) interrupting the analyticity of the series.  A large part of this paper describes a method of determining CLSPs for power expansions of algebraic branches.  The method is based on analytic continuation of the branching sheets over the singular points of $w(z)$. 

 The Resultant of $f(z,w)$ with its derivative $f_w$ is denoted by $R(f,f_w)$.  A point $z_s$ where $R(f,f_w)$ is zero is called a singular point of $w(z)$.  However this does not tell us which branch of $w(z)$ is singular.  In this paper, singular points are sorted by absolute value and labeled $s_i$ with $i$ ranging from one to the total number of finite singular points.  The point at infinity is included if the function is singular there. 

 A finite singularity $(s_i,w_s)$ of $f(z,w)$ is a point in $\mathbb{C}^2$ so may not affect all branch coverings of $w(z)$ over $s_i$  unless $w(z)$ fully-ramifies at $s_i$ into a $n$-cycle branch.  For example, a $10$-degree function at $s_i$ may only ramify into a single $2$-cycle branch and eight analytic $1$-cycles.  In this case, the $2$-cycle covering is singular.  The single-cycle branches are not analytically affected by the singularity.  However, if the function fully-ramifies into a $10$-cycle branch represented by ten series each with power expansions in terms of $z^{1/10}$, then all coverings of this branch would be affected by the singularity.  It is for this reason a power expansion of an algebraic function often has a region of convergence extending beyond the nearest singular point:  some of the branch coverings may simply not be affected by the singularity.  Only when the branch covering, and by association the corresponding power expansions, encounters a convergence-limiting singular point does the analytic domain of the branch and convergence region of its power expansions become established.  
 
 The procedure for finding a convergence-limiting singular point $s_l$ for an expansion around $s_b$ is to analytically extend branch sheets across successively distant singular points over analytic $1$-cycle branches until the continuation reaches either a multi-cycled branch sheet or a single cycle sheet with a pole.

\subsection{Necessary and sufficient conditions for continuing a branch across a singular point:}
\label{sub:section0041}
\begin{enumerate}
	\item In order to analytically continue an $n$-cycle branch from one singular point to the next nearest singular point, the next singular point must have at least $n$ single-cycle analytic branches to support continuity, i.e., $1$-cycle branches which do not have poles.
	\item All branch sheets of an $n$-cycle branch must continue onto analytic $1$-cycle branches to be continuous over a singular point. 
	\item   Branch sheets of an $n$-cycle branch may continue across different singular points but the analytic region of the branch is established upon encountering the nearest multi-cycle branch sheet or branch sheet with a pole.  The first singular point in which this occurs is the CLSP for the associated set of conjugate Puiseux series and establishes their radii of convergence. 
	\label{enum:enum1}
\end{enumerate}

\section{Procedure}
 
The following steps are used to identify convergence-limiting singular points:
	
	\begin{enumerate}
	\item Compute the set $\{s_i\}$ of singular points of $w(z)$.
		\item Compute Puiseux series centered at each singular point to a desired working precision and identify branch cycle types according to the categories in Section (\ref{sec:bt001}).
		\item Construct analytically-continuous routes between a base singular point,$s_b$, and remaining singular points, $s_i$,  in order of increasing distance from the base singular point according to Figure \ref{figure:figure1}.
		\item  Compute the value of each branch sheet at points A and D in Figure \ref{figure:figure1} to a desired precision using the Puiseux series computed above.
		\item Identify removable singular points and adjust the route according to the second diagram in Figure \ref{figure:figure1}. 
		\item Numerically integrate $w(z)$ from A to D shown in Figure \ref{figure:figure1}.
		\item Identify branch sheets over the base singular point, $s_b$ that are analytically-continuous across successive singular points until all branch sheets have encountered convergence-limiting singular points.
		\end{enumerate}

\subsection{Computing singular points}

Singular points of $w(z)$ are easily computed to a desired precision by solving for the zeros of the resultant $R(f,f_w)$ using Mathematica's $\texttt{NSolve}$.  Each Puiseux series will have a convergence-limiting singular point,$s_l$, which determines the radius of convergence of the series and for numerical work, singular points were computed with $400$ to $5000$ digits of precision depending on the test case.    Power series with negative exponents represent branches with poles.  Poles of $w(z)$ are the roots of $a_n(z)$.  The singular points are ordered with increasing distance from a selected base singular point $s_b$.  

\subsection{Computing Puiseux series around singular points}

Mathematica 12.1 and later versions include $\texttt{AsymptoticSolve}$ for computing Puiseux series of algebraic curves defined by (\ref{eqn:eqn001}).  The function attempts to compute the series exactly in terms of rational coefficients when possible.  Otherwise the expansions are done numerically and the precision of the series is limited by the precision of the singular points.  In some cases in this study, $\texttt{AsymptoticSolve}$ was not able in a reasonable time to compute a set of Puiseux series with a sufficient number of terms.  In those cases, custom code written by this author was used to compute the Puiseux series.

Once the Puiseux series have been computed, the next step is to determine the cycle size of each series.  This is determined by taking the least common denominator (LCD) of the powers of $z$ for each series.  The conjugate series associated with these branches are easily identified by computing the LCD of all exponents in a series.  Series having the same LCD belong to the same type of conjugate class.  If there are more than $k$ series with LCD $k$, then each series can be conjugated to determine which series belongs to each conjugate-$k$ class. 
\subsection{Constructing an analytically-continuous route between deleted domains of singular points}
\label{sec:sec005}
Since the Puiseux series for (1) converge at least up to the nearest singular point, an analytically-continuous route can be created between a deleted neighborhood of one singular point and another singular point.  The routes are created for successively distant singular points relative to a base singular point.  In the method employed below, only branch sheets that are continuous across previous singular points are extended across the next successive singular point.  The path of continuation is shown in the first diagram of Figure \ref{figure:figure1} from point $A$ to $D$.  For each singular point, the associated Puiseux series converges in a domain at least equal to the nearest singular point.  In order to obtain precise results at point $A$ and $D$ using the Puiseux series, a deleted domain with radius of $1/3$ the distance to the nearest singular point is chosen and shown as the circles around each singular point in the diagram.  A straight line path between the two domains is computed according to (\ref{eqn:eqn1000}) and then each branch sheet of $s_b$ is  numerically integrated over the path from point $A$ to point $D$ onto a branch sheet of $s_n$.  Each $s_b$ branch sheet is then checked for analytic continuity as per Section \ref{sub:section0041}.  

After $s_n$ has been checked this way and branch sheets have been found to be continuous, the next sequential singular point is tested with the continuous branch sheets.  However there is the possibility of attempting to continue a branch sheet to another singular point when a removable singular point is in the path of integration.  Numerical integration will fail over a removable singularity even though the function is analytic because Equation (\ref{eqn:eqn002}) is used for the derivative and at a removable singular point, this quotient is indeterminate.   In this case, the integration path is split into paths $\beta_1$ and $\beta_2$ shown in the second diagram of Figure \ref{figure:figure1}.
\begin{figure}[H] 
\centering
\includegraphics[scale=0.75]{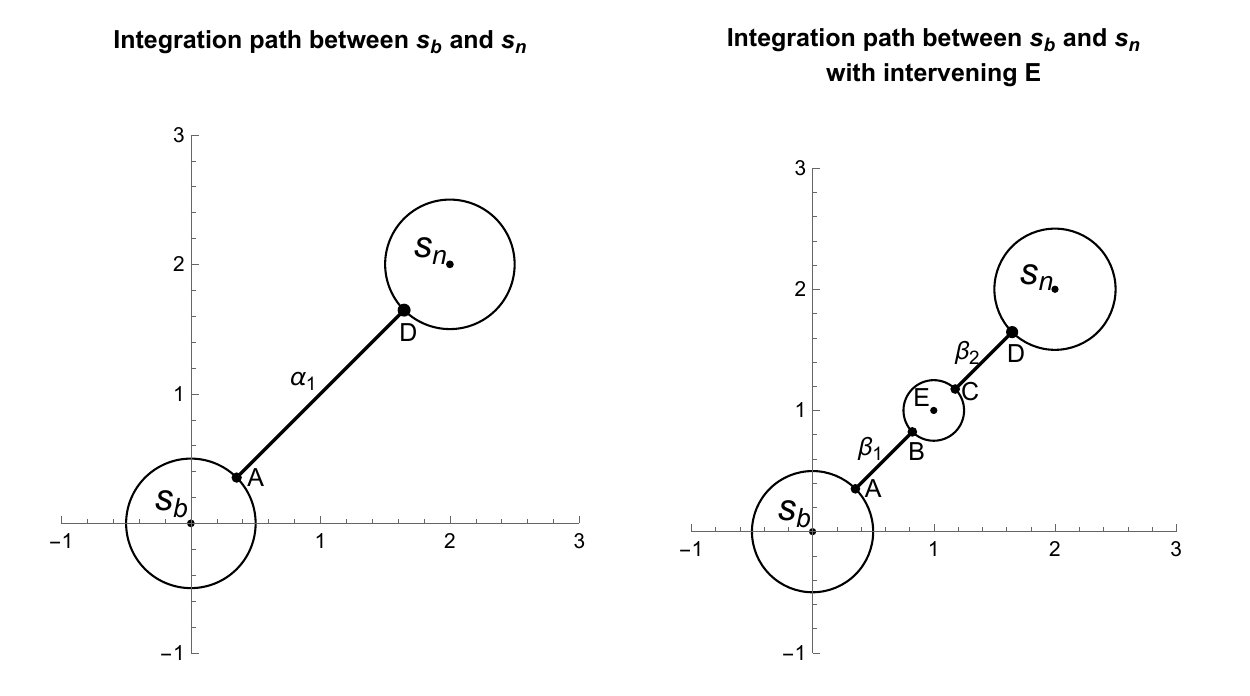}
\caption{Continuation Path}
\label{figure:figure1}
\end{figure}

\begin{equation}
z(t)=A(1-t)+D t;\quad 0\leq t\leq 1
\label{eqn:eqn1000}
\end{equation}

\subsection{Integrating the monodromy differential equation}

For each singular point, Puiseux series of  $128$ to $2040$ terms were computed.  Since the series are convergent in the deleted neighborhoods, the terms of the series eventually become and remain smaller than a minimum desired accuracy.  For this study, at least the number of terms needed to obtain an accuracy of $10^{-10}$ was used.  Often, this was less than $100$ terms.  

Each branch sheet is analytically-continued from the $s_b$ domain to the the next $s_n$ domain according to Figure \ref{figure:figure1} via numerical integration of the monodromy differential equation:

Given $f(z,w)=a_0(z)+a_1(z)w+\cdots+a_n(z)w^n=0$, then $\displaystyle \frac{dw}{dt}=-\frac{f_z}{f_w}\frac{dz}{dt}$.  A set of initial value problems is next set up, one for each branch sheet with initial values at $A$ for each branch sheet and then integrated along the path in Figure \ref{figure:figure1} to the point $D$.  For example, the following $n$ initial value problems for a $n$-degree function are initially created for a route $\alpha_1$ integration path:
$$
\frac{dw}{dt}=-\frac{f_z}{f_w}\frac{dz}{dt}, \quad\; w_i(0)=P_i(A-s_b); \quad i=1,2,\cdots,n 
$$
with  $z(t)=A(1-t)+D t;\quad 0\leq t\leq 1$ and each $P_i(z)$ is a Puiseux series at $s_b$.

 Once precise values of the function were determined at the integration end points, an association between the base branch sheets and the branch sheets about $s_n$ could be made according the the requirements listed in Subsection \ref{sub:section0041} in order to determine which base sheets were analytically continuous over the current singular point.    

\subsection{Illustrative example of continuation across multiple singular points}

\begin{figure}[H]
\centering
\includegraphics[scale=0.7]{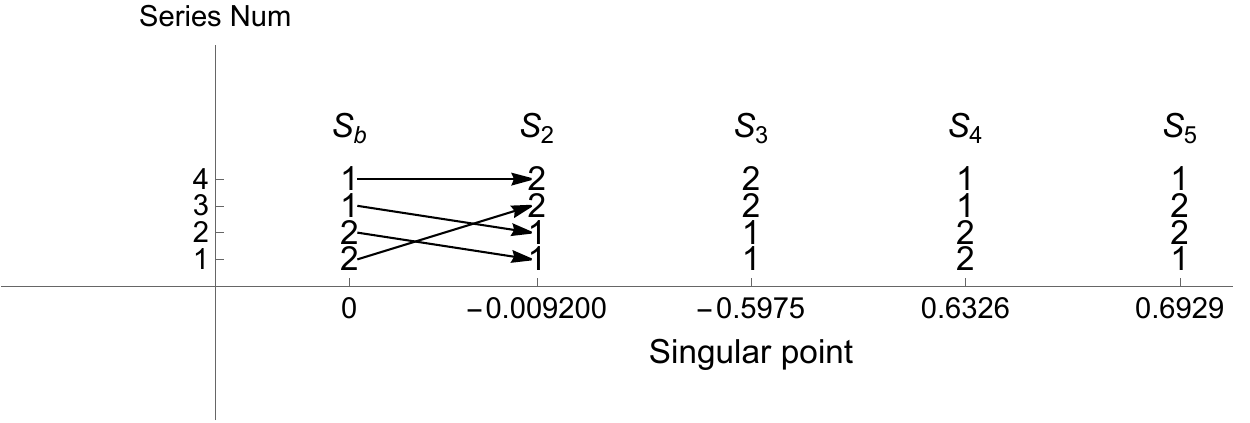}
\caption{Branch links between $s_b$ and $s_2$}
\label{figure:figure2}
\end{figure}

 This section describes the analytic-continuation process graphically using the degree-$4$ function  of  Test Case 2 below.  This function has  four coverings over the complex $z$-plane in deleted neighborhoods of singular points.  The singular point at the origin, $s_b$, ramifies into a $2$-cycle branch and two $1$-cycle branches with each branch sheet represented by a conjugate Puiseux series.  The notation $\{1,2,2,1\}$ refers to the associated sorted Puiseux series of each branch sheet in the series expansion set $\{P_i\}$ at $s_b$:  the first Puiseux series is an analytic $1$-cycle branch (the first number in the set), the second series is a sheet of the $2$-cycle branch, the third is the second sheet of the $2$-cycle branch, and the last Puiseux series is an analytic $1$-cycle.  The nearest singular point in this example is $s_2=-0.0092$ with branching cycles of $\{1,1,2,2\}$ with the $1$-cycles analytic (without a pole).  Figure \ref{figure:figure2} shows the branch sheet types above the base singular point and four nearest singular points.  The vertical axis is the Puiseux series number of each branch sheet. 

 Starting from series $1$ of $s_b$, the figure shows one $2$-cycle sheet of $s_b$ continuing onto a $2$-cycle sheet of $s_2$ (index number $3$ in the diagram).  Series $2$ of $s_b$, the second sheet of the $2$-cycle, continues onto a $1$-cycle sheet of $s_2$.  The third sheet is a $1$-cycle continuing onto a $1$-cycle of $s_2$, and the fourth series of $s_b$ is a $1$-cycle continuing onto a $2$-cycle sheet of $s_2$.  Thus, only the second and third sheets of $s_b$ can continue across $s_2$ over $1$-cycle sheets.

\begin{figure}[H]
\centering
\includegraphics[scale=0.7]{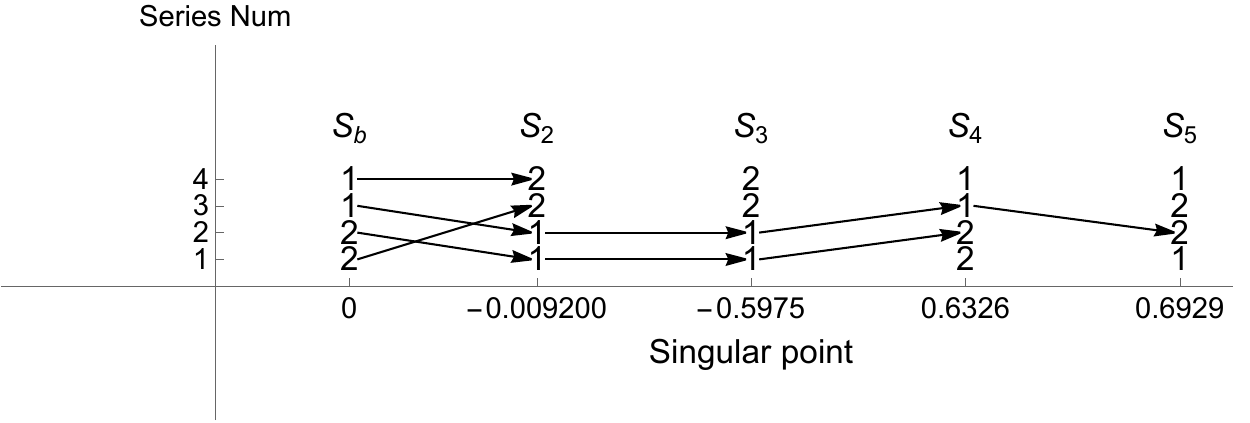}
\caption{Complete continuation branching over $s_b$}
\label{figure:figure4}
\end{figure}

 We next attempt to continue sheets $2$ and $3$ across $s_3$.  This is shown in Figure \ref{figure:figure4}.  Sheets $2$ and $3$ of $s_b$ continue onto $1$-cycle branches of $s_3$ as shown by the arrows between $s_2$ and $s_3$.  At $s_4$, the continuation of sheet $2$ of $s_b$ stops as it continues onto a $2$-cycle sheet.  Sheet $3$ however continues onto another $1$-cycle sheet so that we move over the next nearest singular point $s_5$ and finally, by following the arrows, sheet $3$ of $s_b$ terminates onto a $2$-cycle sheet of $s_5$. 
 
 We can now identify the convergence-limiting singular points of each of the Puiseux series of $s_b$.  However, we need to be careful about multiple-cycle branch sheets:  In Figure \ref{figure:figure4} individual sheets of the $2$-cycle branch of $s_b$ reached different convergent-limiting singular points.  Series $1$ terminated at $s_2$ and series $2$ terminated at $s_4$.  However, both series are in the same conjugate class so the absolute value of conjugate terms in each series are equal.  By the Root Test, both series will converge to the same value which must be the minimum of $|s_2|$ and $|s_4|$.  Therefore, the (exact) radius of convergence of the $2$-cycle series of $s_b$ is $|s_2|$ although numerically, the valued is limited by the precision of the singular point.  From the figure, the convergence radius of series $4$ is also $|s_2|$, and the convergence radius of series $3$ is $|s_5|$.  

\subsection{Applying the Root Test to check the Results}

  The Root Test can be used to approximate the radius of convergence of Puiseux series and used as a check of the analytic continuation results of the previous section.  To use the Root Test, the standard definition is modified to include the branch cycle size:
\begin{equation}
R=\frac{1}{\displaystyle \liminf_{k\to \infty}\; |a_k|^{\frac{c}{m_k}}}
\label{eqn:eqnss2}
\end{equation}
 
where $c$ is the cycle size of the series, and the set $\{m_i\}$ is the set of exponent numerators under a least common denominator.  For example, the terms
$$
a_0+a_1 z+a_2 z^{3/2}+a_3 z^{9/4}+z^{3}
$$

would have the set $\{m_i\}=\{0,4,6,9,12\}$.  Then the $\liminf$ of each Puiseux series can be approximated by forming the set:

$$
S=\bigg\{\left(\frac{1}{m_k},\frac{1}{|a_k|^{\frac{c}{m_k}}}\right)\bigg\}
$$
and numerically extrapolating $\displaystyle \lim_{k\to\infty} S$ using a sufficient number of trailing points of $S$.
For example, consider $512$ terms of a $2$-cycle series:

$$
a_1 z^{1/2}+a_2 z^{2}+a_3 z^{5/2}+a_4 z^{3}+\cdots+a_{512}z^{512}
$$
Therefore we have the sequence:
$$
S=\biggr\{\left(\frac{1}{1},\frac{1}{|a_1|^{2}}\right),\left(\frac{1}{4},\frac{1}{|a_2|^{1/2}}\right),\left(\frac{1}{5},\frac{1}{|a_3|^{2/5}}\right),\left(\frac{1}{6},\frac{1}{|a_3|^{1/3}}\right),\cdots,\left(\frac{1}{1024},\frac{1}{|a_{512}|^{1/512}}\right)\biggr\}
$$
The next step is to approximate $R$ by fitting a suitable curve to the trailing terms of $S$.   In some cases, a large number of terms are needed before the sequences settles into a trend.  Up to $2040$ terms were used in Test Case $2$ and only after several hundred terms did the sequences start to become monotonic towards a limiting point.  Linear or quadratic fits were used to fit the data and extrapolated to $\displaystyle \lim_{k\to\infty} S$.  The value computed by the Root Test was then compared to the value determined by analytic continuation of branch sheets.  In all cases studied, the percent error between the Root Test and analytic continuation was no more than $1.3\%$ and likely could be lowered by using more series terms.  An example plot showing this analysis is shown in Figure \ref{figure:figure8}.
%
%

\section{Test Case 1}

 Singular points for Equation (\ref{eqn:testcase1eqn001}) were computed to $400$ digits of precision via $\texttt{NSolve}$,  and expansions around all singular points were computed to $200$ digits of precision. The expansions around each singular point were generated with $512$ terms and the branching analytically continued as per Section 5 until reaching CLSPs.  Between $7$ and $42$ terms of the series were needed to analyze the branches at the markers $A$ and $D$ in Figure \ref{figure:figure1} to at least $10$ digits of accuracy.  However, $512$ terms were used in the Root Test to better approximate the convergence radii as a check of CLSP results computed by analytic continuation.

\begin{equation}
f_1(z,w)=-z^3 + (z + z^2 + z^3) w + 2 z^2 w^2 + (-1 + z + z^3) w^3=0
\label{eqn:testcase1eqn001}
\end{equation}

Puiseux series up  to degree $512$ and center at the origin can be generated for this function with the following Mathematica code:
\begin{center}
\begin{minipage}{0.85\linewidth}
\begin{lstlisting}[style=myMathematica,
frame=single,
caption=Code to compute Puiseux series of \ref{eqn:testcase1eqn001} at zero,
label=code3]
theFunction = -z^3 + (z + z^2 + z^3) w + 2 z^2 w^2 + (-1 + z + z^3) w^3;
pSeries = w /. AsymptoticSolve[theFunction == 0, w, {z, 0, 512}];
	\end{lstlisting}
	\end{minipage}
	\end{center}

This produces the following terms for each power series:
\begin{equation}
\begin{aligned}
P_1(z)&=-1. z^{1/2}-1. z^{3/2}+0.5 z^2-1. z^{5/2}+1.5 z^3-1.625 z^{7/2}+z^4-2.625 z^{9/2}+\cdots+a_1 z^{512}\\
P_2(z)&=z^2-1. z^3+4. z^7-9. z^8+9. z^9-7. z^{10}-12. z^{11}+91. z^{12}-222. z^{13}+337. z^{14}+\cdots+a_2 z^{512}\\
P_3(z)&=z^{1/2}+z^{3/2}+0.5 z^2+z^{5/2}+1.5 z^3+1.625 z^{7/2}+z^4+2.625 z^{9/2}+\cdots+a_3 z^{512}\\
\label{eqn:case1eqn2}
\end{aligned}
\end{equation}

 Note (\ref{eqn:case1eqn2}) contains two series with $k/2$ powers, and one series with integer powers of $z$.  The two series with fractional exponents represent a conjugate set of power series for a $2$-cycle branch with each series representing a single-valued sheet of the branch.  As an example of series conjugation, consider the first series in (\ref{eqn:case1eqn2}):
$$
P_1(z)=-1. z^{1/2}-1. z^{3/2}+0.5 z^2-1. z^{5/2}+1.5 z^3-1.625 z^{7/2}+z^4-2.625 z^{9/2}+\cdots+a_1 z^{512}.
$$
In order to generate $p_3(z)$ by conjugation, apply (\ref{eqn:eqn301}) to $p_1(z)$:
$$
\begin{aligned}
&(-1)\left(e^{\pi i}\right)^1 z^{1/2}+(-1)\left(e^{\pi i}\right)^3 z^{3/2}+(0.5)\left(e^{\pi i}\right)^4 z^{4/2}\\
&-\left(e^{\pi i}\right)^5 z^{5/2}+1.5\left(e^{\pi i}\right)^6 z^{6/2}-1.625 \left(e^{\pi i}\right)^7 z^{7/2}\\
&+\left(e^{\pi i}\right)^8 z^{8/2}-2.625\left(e^{\pi i}\right)^9 z^{9/2}+\cdots+a_1 \left(e^{\pi i}\right)^{1024} z^{1024/2}\\
&=P_3(z)
\end{aligned}
$$

\begin{table}[ht]
\caption{Singular point branch expansions for (\ref{eqn:testcase1eqn001})} 
\centering
\begin{tabular}{|c|c|l|c|>{$}c<{$}|}
\hline
\text{$s_b$} & \text{$s_b$ Value} & \text{(Branch Type,CLSP)} & \text{$\%$ Error} \\ \hline
$s_{1}$ & $0$ & $(V_{2},s_{2})$, $(E,s_{2})$& $0.65$\\$s_{2}$ & $-0.3582 - 0.2530 i$ & $(V_{2},s_{1})$, $(T,s_{1})$& $0.57$\\$s_{3}$ & $-0.3582 + 0.2530 i$ & $(V_{2},s_{1})$, $(T,s_{1})$& $0.57$\\$s_{4}$ & $0.1961 - 0.5259 i$ &$(T,s_{8})$, $(V_{2},s_{8})$& $0.56$\\$s_{5}$ & $0.1961 + 0.5259 i$ &$(T,s_{9})$, $(V_{2},s_{9})$& $0.56$\\$s_{6}$ & $0.6823$ & $(L^{-1},s_{7})$, $(T,s_{7})$, $(T,s_{4})$& $0.63$\\$s_{7}$ & $0.7492$ & $(V_{2},s_{6})$, $(T,s_{4})$& $0.52$\\$s_{8}$ & $-0.1440 - 0.9401 i$ & $(V_{2},s_{10})$, $(T,s_{4})$& $0.42$\\$s_{9}$ & $-0.1440 + 0.9401 i$ & $(V_{2},s_{11})$, $(T,s_{5})$& $0.42$\\$s_{10}$ & $-0.3412 - 1.1615 i$ & $(L^{-1},s_{8})$, $(T,s_{8})$, $(T,s_{12})$& $0.43$\\$s_{11}$ & $-0.3412 + 1.1615 i$ & $(L^{-1},s_{9})$, $(T,s_{9})$, $(T,s_{13})$& $0.43$\\$s_{12}$ & $-0.843 - 1.560 i$ &$(T,s_{8})$, $(V_{2},s_{10})$& $0.43$\\$s_{13}$ & $-0.843 + 1.560 i$ &$(T,s_{9})$, $(V_{2},s_{11})$& $0.43$\\

\hline
\end{tabular}
	\label{table:tabletest}
\end{table}

 Table \ref{table:tabletest} lists the branch types and CLSP's for singular point power expansions of (\ref{eqn:testcase1eqn001}) centered at each singular point $s_b$.  These CLSP's are relative to the expansion center:  if a branch expansion about $s_n$ has a CLSP of $s_k$ then the radius of convergence for the associated power series is $R=|s_n-s_k|$.  The third column in the table is the percent error between the convergence radius determined by the root test and by CLSPs.  Note in particular there is an $E$ branch type at the origin. 

  First consider the expansion around the origin, $s_1=0$, with branch codes $((V_2,s_2),(E,s_2))$.  This set of Puiseux series at this singular point consists of a $2$-cycle $V_2$ branch with radius of convergence equal to $R=|s_1-s_2|=|s_2|\approx 0.438$, and a $1$-cycle $E$ branch with a removable singularity having radius of convergence also $|s_2|$.	  Figure \ref{figure:case1Figure1} shows the real parts of these branches.  Note the two sheets of $V_2$ over the complex plane and the single sheet of the $E$ branch.  These figures were plotted using the first $50$ terms of Puiseux expansions for the branches with $0<r<0.43$.  See Listing \ref{code4} for example code for plotting branches using their Puiseux series.
 
 Next consider the expansion around $s_5$ with branch expansions $(T,s_9),(V_2,s_9)$.  The first branch is a $1$-cycle Taylor expansion with radius of convergence $R=|s_5-s_9|\approx 0.536$.  The second expansion is a $2$-cycle $V_2$ branch (with vertical tangent at center) having a power expansion with radius of convergence $R=|s_5-s_9|$.
	
\begin{figure}[H] 
\centering
\includegraphics[scale=0.7]{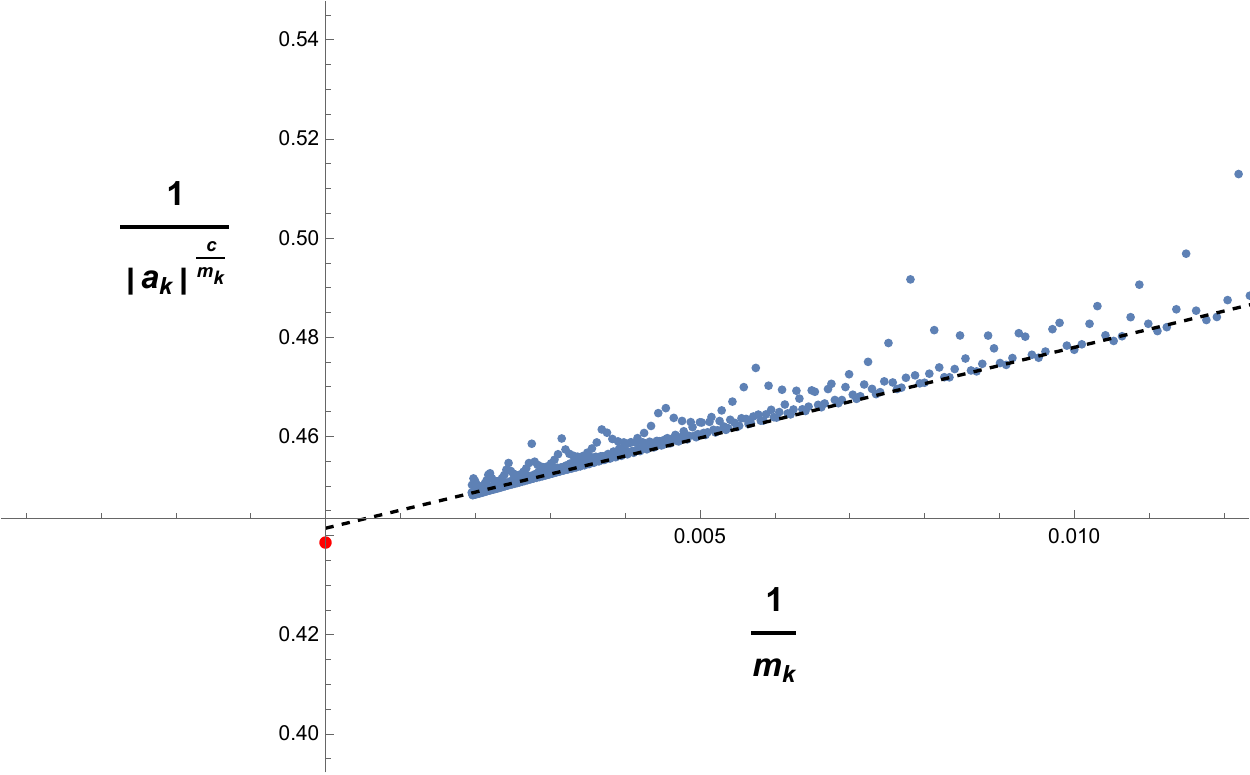}
\caption{Root test for $1$-cycle branch over origin}
\label{figure:figure8}
\end{figure}

\begin{figure}[H]
\centering
\includegraphics[scale=0.7]{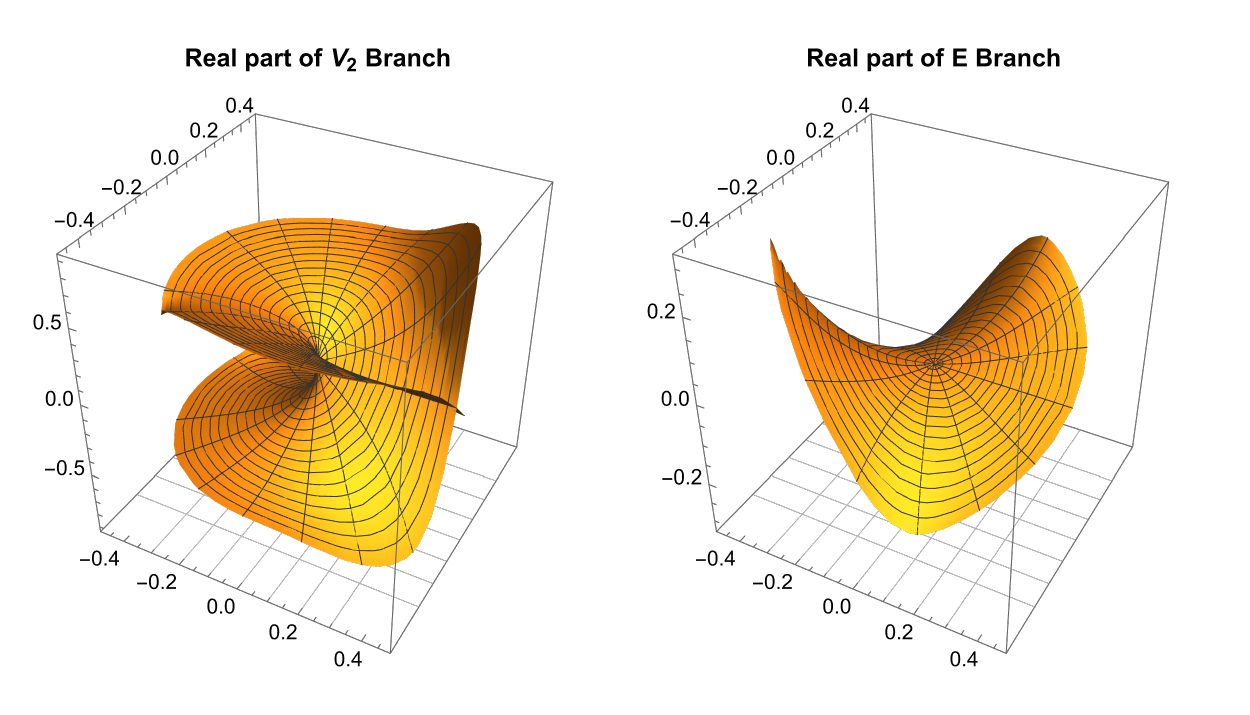}
\caption{Real part of $V_2$ and $E$ Branches at the origin}
\label{figure:case1Figure1}
\end{figure}

 Figure \ref{figure:figure8} is a sample plot using a quadratic fit of Root Test data to the $1$-cycle $(E,s_2)$ series at the origin.  $512$ terms of the series were used to generate this plot. In Figure \ref{figure:figure8}:
\begin{enumerate}
	\item The scatter points represents the set of points in $S$.
	\item The dashed black line is a quadratic fit to the lower boundary of scatter points.  This line is use to approximate $\liminf$.
		\item The intersection of the dashed black line with the vertical axis is an approximation to the radius of convergence for this series.
	\item A red dot on the vertical axis represents the radius of convergence determined by analytic continuation.
	\end{enumerate}

%
%
%

\section{\textbf{Test Case 2}}

\begin{equation}
\begin{aligned} 
f_2(z,w)&=\left(-z^2+z^3\right)\\
&+\left(-4 z+3 z^2\right) w\\
&+\left(-z^3-9 z^4\right) w^2\\
&+\left(-2+8 z+4 z^2-4 z^3\right) w^3\\
&+\left(6-8 z^2+7 z^3+8 z^4\right) w^4=0
\end{aligned}
\label{eqn:eqn005}
\end{equation}

Table \ref{table:table66} list the singular points and branch codes for (\ref{eqn:eqn005}) and includes the point at infinity. The values were computed to $400$ digits of precision using between $128$ and $512$ terms.

\begin{table}[ht]
\caption{Singular points for (\ref{eqn:eqn005})} 
\centering 
\begin{tabular}{|c|c|c|c|>{$}c<{$}|}\hline
\text{$s_b$} & \text{$s_b$ Value} & \text{(Branch Type,CLSP)} & \text{$\%$ Error} \\ \hline
$s_{1}$ & $0$ & $(V_{2},s_{2})$, $(E,s_{5})$, $(E,s_{2})$& $0.80$\\$s_{2}$ & $-0.009200$ &$(T,s_{1})$, $(T,s_{5})$, $(V_{2},s_{1})$& $0.57$\\$s_{3}$ & $-0.5975$ &$(T,s_{16})$, $(T,s_{12})$, $(V_{2},s_{12})$& $0.62$\\$s_{4}$ & $0.6326$ & $(V_{2},s_{5})$, $(T,s_{5})$, $(T,s_{17})$& $0.63$\\$s_{5}$ & $0.6929$ &$(T,s_{4})$, $(T,s_{17})$, $(V_{2},s_{4})$& $0.90$\\$s_{6}$ & $0.6447 - 0.5028 i$ & $(L^{-1},s_{14})$, $(T,s_{14})$, $(T,s_{5})$, $(T,s_{17})$& $0.43$\\$s_{7}$ & $0.6447 + 0.5028 i$ & $(L^{-1},s_{15})$, $(T,s_{15})$, $(T,s_{5})$, $(T,s_{18})$& $0.43$\\$s_{8}$ & $0.2964 - 0.7975 i$ &$(T,s_{14})$, $(T,s_{10})$, $(V_{2},s_{10})$& $0.43$\\$s_{9}$ & $0.2964 + 0.7975 i$ &$(T,s_{15})$, $(T,s_{11})$, $(V_{2},s_{11})$& $0.43$\\$s_{10}$ & $-0.0729 - 0.8528 i$ &$(T,s_{24})$, $(T,s_{8})$, $(V_{2},s_{8})$& $0.57$\\$s_{11}$ & $-0.0729 + 0.8528 i$ &$(T,s_{25})$, $(T,s_{9})$, $(V_{2},s_{9})$& $0.57$\\$s_{12}$ & $-0.8591$ &$(T,s_{16})$, $(T,s_{13})$, $(V_{2},s_{3})$& $0.58$\\$s_{13}$ & $-0.8608$ &$(T,s_{16})$, $(T,s_{12})$, $(T,s_{12})$, $(L^{-1},s_{16})$& $0.44$\\$s_{14}$ & $0.7205 - 0.4925 i$ & $(V_{2},s_{6})$, $(T,s_{5})$, $(T,s_{17})$& $0.42$\\$s_{15}$ & $0.7205 + 0.4925 i$ & $(V_{2},s_{7})$, $(T,s_{5})$, $(T,s_{18})$& $0.42$\\$s_{16}$ & $-0.9016$ & $(V_{2},s_{13})$, $(T,s_{12})$, $(T,s_{12})$& $0.42$\\$s_{17}$ & $0.8593 - 0.4299 i$ &$(T,s_{14})$, $(T,s_{19})$, $(V_{2},s_{14})$& $0.57$\\$s_{18}$ & $0.8593 + 0.4299 i$ &$(T,s_{15})$, $(T,s_{19})$, $(V_{2},s_{15})$& $0.57$\\$s_{19}$ & $0.9666$ &$(T,s_{4})$, $(T,s_{17})$, $(V_{2},s_{5})$& $0.57$\\$s_{20}$ & $-1.1628 - 0.2641 i$ &$(T,s_{22})$, $(T,s_{12})$, $(V_{2},s_{22})$& $1.1$\\$s_{21}$ & $-1.1628 + 0.2641 i$ &$(T,s_{22})$, $(T,s_{12})$, $(V_{2},s_{22})$& $1.1$\\$s_{22}$ & $-1.296$ &$(T,s_{20})$, $(T,s_{21})$, $(V_{2},s_{23})$& $0.43$\\$s_{23}$ & $-1.304$ &$(T,s_{20})$, $(T,s_{21})$, $(T,s_{22})$, $(L^{-1},s_{22})$& $0.44$\\$s_{24}$ & $0.2805 - 1.3743 i$ & $(V_{2},s_{10})$, $(T,s_{8})$, $(T,s_{8})$& $0.57$\\$s_{25}$ & $0.2805 + 1.3743 i$ & $(V_{2},s_{11})$, $(T,s_{9})$, $(T,s_{9})$& $0.57$\\
$s_{\infty}$ & & $(T,s_2),(T,s_5),(V_2,s_2)$ & $1.1$ \\

 \hline
\end{tabular}
	\label{table:table66} 
\end{table}

 The origin contains two $E$ branches which are branches with removable singularities at their centers.  And since there are other singular points on the real axis, these removable points may intersect an integration path as described in Section \ref{sec:sec005}.  Take for example, the expansion around $s_2=-0.0092$.  The sequence of singular points nearest to this point includes $\{0,-0.597,0.63\}$.  Notice if the analytic continuation proceeds to the point $0.63$, the integration path start at $-0.0092$ through the origin, and onto the domain of $0.63$.  In this case, the path is split into the $\beta_1$ and $\beta_2$ paths shown in Figure \ref{figure:figure4}.  

\begin{figure}[H]
\centering
\includegraphics[scale=0.7]{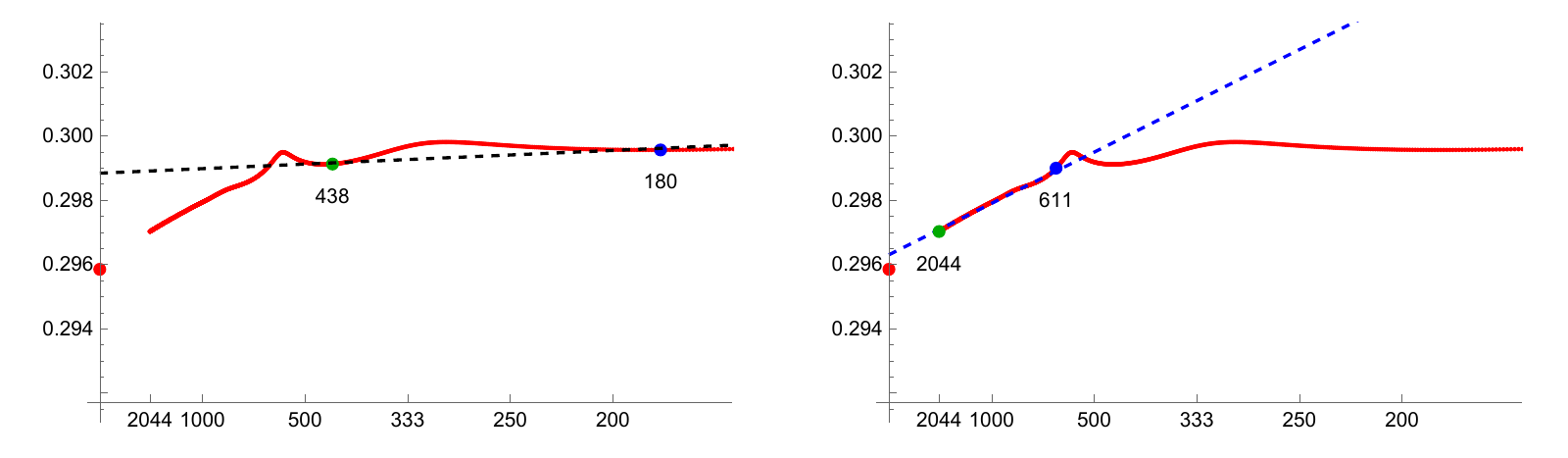}
\caption{Root test for $(T,s_{22})$ branch about $s_{20}$ in Test Case $2$}
\label{figure:figure2001}
\end{figure}

 In Table \ref{table:table66}, the percent error between the Root Test and the method of analytic continuation for the expansions about $s_{20}$, $s_{21}$ and $s_{\infty}$ were reported above $1\%$.  Closer examination explains this unusually high error:  The series for $s_{20}$ were re-generated with $2044$ terms at $5000$ digits of precision. Figure \ref{figure:figure2001} shows the Root Test results using two different sets of results.  The set $\{S\}$ is shown as the red line with the $x$-axis labeled $i$ for the term $m_i$.  On the vertical axis is a red point showing the radius of convergence for this series as determined by analytic continuation.  The first $500$ terms appear to settle into a relatively stable trend.  In the first figure, terms $180$ to $468$ were  used in the Root Test which extrapolated to a value some distance from the value determined by analytic continuation.  However notice the terms make a steep decline at around the $600$'th term.  In the second figure, terms $611$ to $2044$ were fitted to a linear curve resulting in a closer fit:  When $2044$ Puiseux terms were used for the expansions about $s_{20}$, $s_{21}$ and $s_{\infty}$, the Root Test percent error  was less than $0.3\%$.  This figure demonstrates two important results of this study:  
\begin{enumerate}
	\item The sequence $S$ make take some time before settling down into a trend tending to the convergence point.
	\item Using a sufficient number of terms, $2044$ in this case, the Root Test is in close agreement with the results computed by analytic continuation.
	\end{enumerate}

The expansion around infinity was computed by expanding the function
\begin{equation}
g_2(z,w)=z^d f_2(1/z,w)
\label{eqn:testcase2eqn1000}
\end{equation}

around the origin with $d$ the highest power of $z$ in $f_2(z,w)$ so that an expansion of $g(z,w)$ about the origin is an expansion of $f(z,w)$ around infinity.   The branching ramified into two $1$-cycle branches and a $2$-cycle branch:  $(T,s_2),(T,s_5),(V_2,s_2)$.

  However, the branch codes of the expansion at infinity are in reference to the reverse sort order of the singular points of $f(z,w)$ since the singular points are placed in an order of nearest distance from the point at infinity:  An expansion about infinity becomes $s_b=s_1$.  The singular point nearest the point at infinity becomes $s_2$ which is $s_{25}$ relative to the origin.  The fifth nearest singular point to infinity becomes $s_5$ which is $s_{22}$ relative to the origin and so forth.  Therefore, the radius of convergence of the first and third power expansion in Table \ref{table:table66} for the expansions at infinity ($s_{\infty}$) with branch codes $(T,s_2)$ and $(V_2,s_2)$ is $\displaystyle\left|\frac{1}{s_{25}}\right|\approx 0.713$.  And the convergence radius of $(T,s_5)$ taken with respect to infinity is $\displaystyle \left|\frac{1}{s_{22}}\right|\approx 0.713$. 
		
Consider computing the roots of $f_2(10,w)=0$ as a check of the expansions at infinity.  In order to use the expansions at infinity we would substitute $z\to 1/z$ as per the substitution made in (\ref{eqn:testcase2eqn1000}).  We first compute the roots of $f_2(10,w)=0$, and then evaluate each of the four Puiseux series at infinity using $z=1/10$.   Below is a comparison of the computed roots with the series results:
	
	\begin{equation}
	\begin{array}{cc}
	\text{NSolve Results: } & {-1.00362, -0.0986669, 0.101185, 1.04195}\\
	\text{Series Results: } & {-1.00362, -0.0986669, 0.101185, 1.04195}
	\end{array}
	\end{equation}
%
%
%

\section{\textbf{Test Case 3}}
\label{testcase3}
Puiseux expansions  of (\ref{eqn:case3eqn1}) at all singular points were performed with $512$ terms at $800$ to $1000$ digits of precision.  Table \ref{table:case2table001} list the branching parameters.   

\begin{equation}
f_3(z,w)=-1/2+z^2 w+(-2 z^2-z^3) w^2+(-z/2+2 z^2)w^3+(-z)w^4=0
\label{eqn:case3eqn1}
\end{equation}

		\begin{table}[H]
\caption{Singular points for (\ref{eqn:case3eqn1})} 
\centering 
\begin{tabular}{|c|c|c|c|>{$}c<{$}|}\hline
\text{$s_b$} & \text{$s_b$ Value}&\text{(Branch Type,CLSP)} & \text{$\%$ Error} \\ \hline
$s_{1}$ & $0$ & $(P_{4}^{-1},s_{2})$& $0.69$\\$s_{2}$ & $-0.1796 - 0.3499 i$ & $(V_{2},s_{1})$, $(T,s_{1})$, $(T,s_{1})$& $0.30$\\$s_{3}$ & $-0.1796 + 0.3499 i$ & $(V_{2},s_{1})$, $(T,s_{1})$, $(T,s_{1})$& $0.30$\\$s_{4}$ & $0.6957 - 0.0167 i$ &$(T,s_{5})$, $(T,s_{5})$, $(V_{2},s_{1})$& $0.41$\\$s_{5}$ & $0.6957 + 0.0167 i$ &$(T,s_{4})$, $(T,s_{4})$, $(V_{2},s_{1})$& $0.41$\\$s_{6}$ & $-1.1392 - 0.5535 i$ &$(T,s_{2})$, $(T,s_{7})$, $(V_{2},s_{2})$& $0.58$\\$s_{7}$ & $-1.1392 + 0.5535 i$ &$(T,s_{3})$, $(T,s_{6})$, $(V_{2},s_{3})$& $0.58$\\$s_{8}$ & $0.6753 - 1.1369 i$ &$(T,s_{4})$, $(T,s_{5})$, $(V_{2},s_{4})$& $0.56$\\$s_{9}$ & $0.6753 + 1.1369 i$ &$(T,s_{5})$, $(T,s_{4})$, $(V_{2},s_{5})$& $0.56$\\$s_{10}$ & $3.152$ & $(V_{2},s_{4})$, $(T,s_{4})$, $(T,s_{5})$& $0.55$\\$s_{11}$ & $-4.832$ &$(T,s_{2})$, $(T,s_{6})$, $(V_{2},s_{6})$& $0.59$\\
$s_{\infty}$ & & $(E,s_3),(E,s_2),(P_2^{-2},s_2)$ &$0.8$\\
 \hline
\end{tabular}
	\label{table:case2table001} 
\end{table}

 The expansion about the origin fully-ramifies into a $4$-cycle pole of order one, $P_4^{-1}$, with radius of convergence of the associated Puiseux expansions $R=|s_2|\approx 0.3933$.  Once $R$ has been determined, it's a simple matter to plot the real or imaginary sheets of the branch.  Listing \ref{code4} first computes Puiseux series up to order $25$ for the four sheets of the branch.  Since the branch has a pole at its center, the plot is generated in an annulus $0.03933\leq r\leq 0.3933$.  The remaining code generates a table of four plots, one plot for each series in this $4$-cycle conjugate class each representing a single-valued sheet of the branch.  

\begin{center}
\begin{minipage}{0.85\linewidth}
\begin{lstlisting}[style=myMathematica,
frame=single,
caption=Code to plot $\text{Im}\left(P_4^{-1}\right)$ branch at the origin,
label=code4]
theFunction = -1/2+z^2 w+(-2 z^2-z^3)w^2+(-(z/2)+2 z^2)w^3+(-z)w^4
baseSingularPoint=0;
pSeries = w /. AsymptoticSolve[theFunction == 0, w, {z, 0, 25}];
rEnd=0.3933;
rStart=0.1 0.3933;
(* code for face-grids and plot label omitted for brevity *)
 thePlot = Table[
   ParametricPlot3D[{Re[z] + Re[baseSingularPoint], 
      Im[z] + Im[baseSingularPoint], Im[pSeries[[i]]]} /. 
     z -> r Exp[I t], {r, rStart, rEnd}, {t, -\[Pi], \[Pi]}, 
    BoxRatios -> {1, 1, 1}],
   {i, 1, 4}];
 	Show[thePlot,PlotRange->All]
	\end{lstlisting}
	\end{minipage}
	\end{center}
	
	\begin{figure}[H] 
\centering
\includegraphics[scale=0.7]{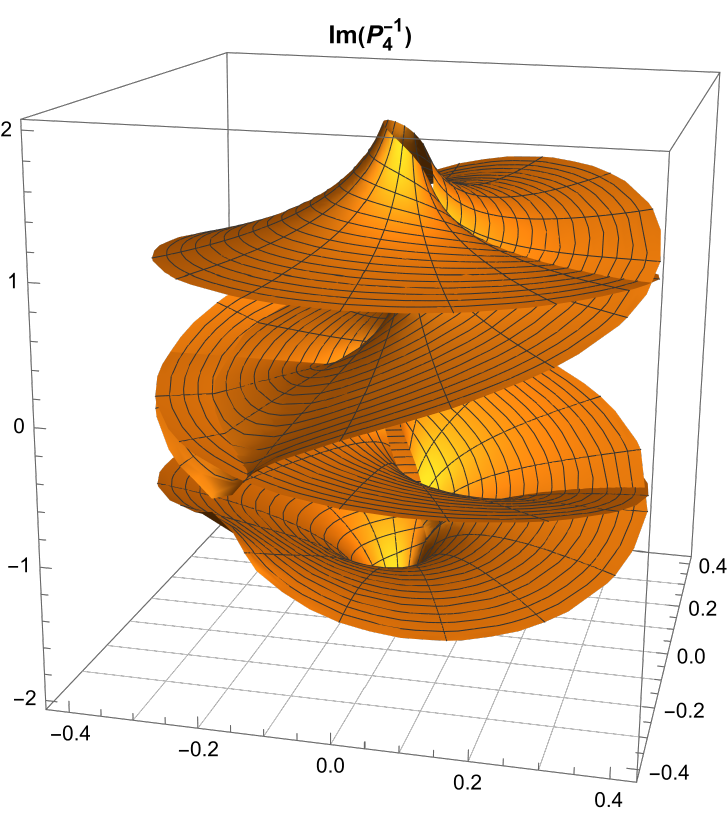}
\caption{Imaginary part of $P_4^{-1}$ branch}
\label{figure:figure3001}
\end{figure}

 Note in Figure \ref{figure:figure3001} the branch surface twisting back on itself after an $8\pi$ route around the $z$-axis giving rise to a four cycle branch consisting of four coverings over the complex $z$-plane.

 And the expansion at infinity has two removable singular points.  We can confirm the removable singular points by analyzing the function
	
	\begin{equation}
	g_3(z,w)=z^3 f_3(1/z,w)=-\frac{z^3}{2}+z w+(-1-2z) w^2+(2z-\frac{z^2}{2})w^3+(-z^2)w^4
	\label{eqn:case2eqn2}
	\end{equation}
		
	and the two partials
	$$
	\begin{aligned}
	\frac{\partial g_3}{\partial z}&=-\frac{3 z^2}{2}+w-2 w^2+(2-z)w^3+(-2 z)w^4 \\
	\frac{\partial g_3}{\partial w}&=z+(-2-4 z)w+(6 z-(3 z^2)/2)w^2+(-4 z^2)w^3
	\end{aligned}
	$$
	
	which are both zero at $\{0,0\}$.
	Next consider the expansions at infinity:
	$$
 \begin{aligned}
P_1(z)&=0.5 z^2+0.25 z^3+0.75 z^4+1.3125 z^5+3.0625 z^6+\cdots\\
P_2(z)&=z-2.5 z^2+5.75 z^3-22.25 z^4+76.9375 z^5+\cdots\\
P_3(z)&=-0.25+1/z-\frac{(0. +1.58114 I)}{\sqrt{z}}+(0. +0.335992 I) \sqrt{z}-0.5 z+\cdots\\
P_4(z)&=-0.25+1/z+\frac{(0. +1.58114 I)}{\sqrt{z}}-(0. +0.335992 I) \sqrt{z}-0.5 z+\cdots
\end{aligned}
$$
and note $\displaystyle \frac{d P_1}{dz}(0)=0$ and $\displaystyle \frac{d P_2}{dz}(0)=1$ so that the point at infinity is a removable singular point over $P_1(z)$ and $P_2(z)$, the two $E$ branches.
%
%
%
\section{Test case 4:}

Expression (\ref{eqn:case4eqn001}) has $180$ singular points, $9$-cycle branches at $\pm i$, fifteen second-order poles at infinity and  ramifies into $1$,$2$,$3$,$4$, and $5$ cycle branches at the origin. 

\begin{equation}
\begin{aligned}
f_4(z,w)&=\left(z^{30}+z^{32}\right)+\left(z^{14}+z^{20}\right) w^5\\
&+\left(z^5+z^9\right) w^9+\left(z+z^3\right) w^{12}+6 w^{14}+\left(2+z^2\right) w^{15}=0
\label{eqn:case4eqn001}
\end{aligned}
\end{equation}

 $512$ terms were generated for each series at the origin with a precision of $1980$.  $128$ terms were generated for the remaining series at a precision of $160$ digits.

\begin{table}[H]
\caption{Branching at the origin and convergence results for (\ref{eqn:case4eqn001})} 
\centering 
\begin{tabular}{|c|c|c|c|c|}
\hline
  \text{Cycle Order}& \text{(Branch Type,CLSP)} & \text{Root Test} & \text{AC} & \text{$\%$ error}\\
	\hline
	1 &$(T,s_{118})$ & $1.099$ & $1.094$ & $0.47$ \\
	2 & $(V_2,s_2)$ & $0.1677$ & $0.1668$ & $0.59$ \\
	3 & $(F_3^4,s_2)$ & $0.168$ & $0.167$ & $0.87$ \\
	4 & $(F_4^9,s_7)$ & $0.510$ & $0.505$ & $1.01$ \\
	5 & $(F_5^{16},s_{27})$ & $0.6488$ & $0.6413$ & $1.2$\\
\hline
 \end{tabular}
	\label{table:table53} 
\end{table}

The branches at the origin were then analytically continued until CLSPs were encountered.  Nearest CLSPs for each conjugate class were then selected and radii of convergences $R=|s_b-s_l|$ were compared to quadratic fits of Root Test results for each series.  Table \ref{table:table53} summarizes the results of this test.  Note in particular the $1$-cycle at the origin continued across $116$ singular points before encountering a CLSP.  Convergence results for the remaining finite singular points can be found at the author's web site:  \href{https://www.jujusdiaries.com/p/radius-of-convergence-of-power.html}{Radius of Convergence of Puiseux series of Algebraic Functions,Part II}

\begin{table}
\caption{Convergence Results for (\ref{eqn:case4eqn001}) at infinity} 
\centering 
\begin{tabular}{|c|c|c|c|c|c|c|}
\hline
   \text{Index} & \text{Size} & \text{Type} & \text{CLSP} & \text{Root Test} & \text{AC} & \text{$\%$ error}\\ \hline
$1$ & $1$ & $L^{-2}$ & $S_{6}$ & $0.551442$ & $0.54979$ & $0.30$\\$2$ & $1$ & $L^{-2}$ & $S_{3}$ & $0.540964$ & $0.538936$ & $0.38$\\$3$ & $1$ & $L^{-2}$ & $S_{2}$ & $0.540964$ & $0.538936$ & $0.38$\\$4$ & $1$ & $L^{-2}$ & $S_{3}$ & $0.541024$ & $0.538936$ & $0.39$\\$5$ & $1$ & $L^{-2}$ & $S_{2}$ & $0.541024$ & $0.538936$ & $0.39$\\$6$ & $1$ & $L^{-2}$ & $S_{11}$ & $0.561659$ & $0.559514$ & $0.38$\\$7$ & $1$ & $L^{-2}$ & $S_{10}$ & $0.561659$ & $0.559514$ & $0.38$\\$8$ & $1$ & $L^{-2}$ & $S_{14}$ & $0.573164$ & $0.571386$ & $0.31$\\$9$ & $1$ & $L^{-2}$ & $S_{15}$ & $0.573164$ & $0.571386$ & $0.31$\\$10$ & $1$ & $L^{-2}$ & $S_{16}$ & $0.581556$ & $0.579324$ & $0.39$\\$11$ & $1$ & $L^{-2}$ & $S_{17}$ & $0.581556$ & $0.579324$ & $0.39$\\$12$ & $1$ & $L^{-2}$ & $S_{22}$ & $0.628151$ & $0.622668$ & $0.88$\\$13$ & $1$ & $L^{-2}$ & $S_{23}$ & $0.628151$ & $0.622668$ & $0.88$\\$14$ & $1$ & $L^{-2}$ & $S_{26}$ & $0.6324$ & $0.629759$ & $0.42$\\$15$ & $1$ & $L^{-2}$ & $S_{27}$ & $0.6324$ & $0.629759$ & $0.42$\\
\hline
 \end{tabular}
	\label{table:table4001} 
\end{table}

\begin{figure}[H]
\centering
\includegraphics[scale=0.7]{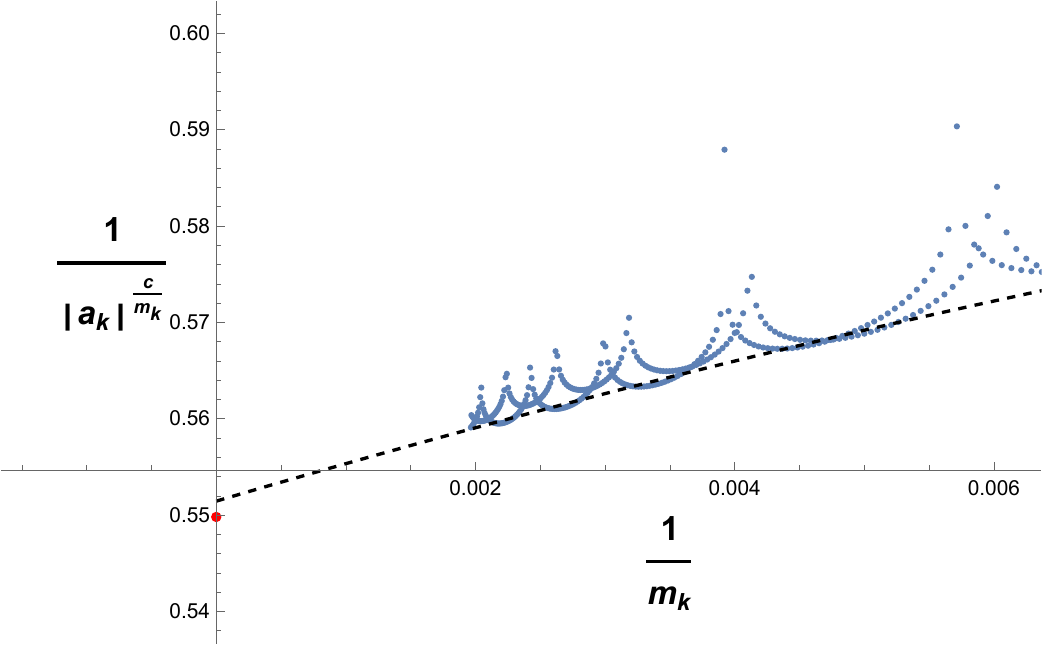}
\caption{Root test for first $L^{-2}$ series in Table \ref{table:table4001}}
\label{figure:figure4001}
\end{figure}

 In order to expand $w_4(z)$ at infinity, the function $g_4(z,w)=z^{32}f_4(1/z,w)$ is expanded at the origin.  $g_4$ ramified into fifteen $1$-cycle poles of order $2$.  Power series up to order $512$ were generated and radii of convergences were computed using both the analytic continuation method and the Root Test.   Convergence results are shown in Table \ref{table:table4001}.  Figure \ref{figure:figure4001} is a Root Test fit of the first $L^{-2}$ series in Table \ref{table:table4001} showing a quadratic fit to the lower limit of points and a red point on the vertical axis showing the radius of convergence for this series computed by analytic continuation.

In order to numerically check the expansions at infinity, choose a value of $z$ outside the domain of finite singular points, say $z_1=20+25i$.  Next compute the values of $f_4(z_1,w)=0$ via $\texttt{NSolve}$.  Then compute $w_i(1/z_1)$ using the power expansions of $f(z,w)$ at infinity.  The $512$ terms of the $S_{\infty}$ series at $1/z_1$ were computed to $1980$ digits of precision and agreed with the roots  to $635$ digits. 

\section{Link to further information}
  Readers interested in further pursuing this subject can find more information at the author's website:  \href{https://www.jujusdiaries.com/p/computing-radii-of-convergence-of-power.html}{Radius of Convergence of Puiseux series of Algebraic Functions}.

\section{Conclusions}
\begin{enumerate}
	\item One might wonder why not just use the Root Test to determine the radius of convergence of the series.  The reason for using analytic continuation is that the deleted domains around singular points may have radii very close to one another requiring an impractical number of terms to resolve the difference between two radii.  And Test Case 2 demonstrated the difficulty of determining an appropriate range of terms to use in the Root Test. The method of analytic continuation precisely identifies the CLSP of a series.  The Root Test was only used as a check of the results.
	\item Since the method relies on numerical integration to analytically continue a branch from one singular point to another, there is the possibility of numerical error causing incorrect results.  However the integration paths are chosen to be mostly smooth although there is the possibility of the integration path traveling very close to a singular point causing a wide oscillation of the numerical results leading to an incorrect CLSP.  Further improvements of the software in this area could identify such cases and minimize errors.    
	\item The method is limited by the ability to compute Puiseux series and successfully integrate the function along the continuation paths.  Functions of degree $15$ and lower were used as test cases.  Further work with more complex functions and benchmarks on computations would be another area for improvement.
	\item The method presented here is a simple solution for determining radii of convergence of power expansions about singular points of the functions studied in this paper.  Puiseux expansions are relatively easily generated for low-degree functions, and once convergence radii are determined, the power expansions become a more robust tool for further exploring these functions.
	
	\end{enumerate}

\pagebreak

\appendix

\section{}
\label{appendix:appA}
\begin{center}
Example real surfaces of the six branch types studied in this paper
\end{center}

 In the figures below, types $P$ and $L$ branches extend to infinity at their centers and so the plots for those branches were clipped for display purposes.
\begin{figure}[H] 
\centering
\includegraphics[scale=0.75]{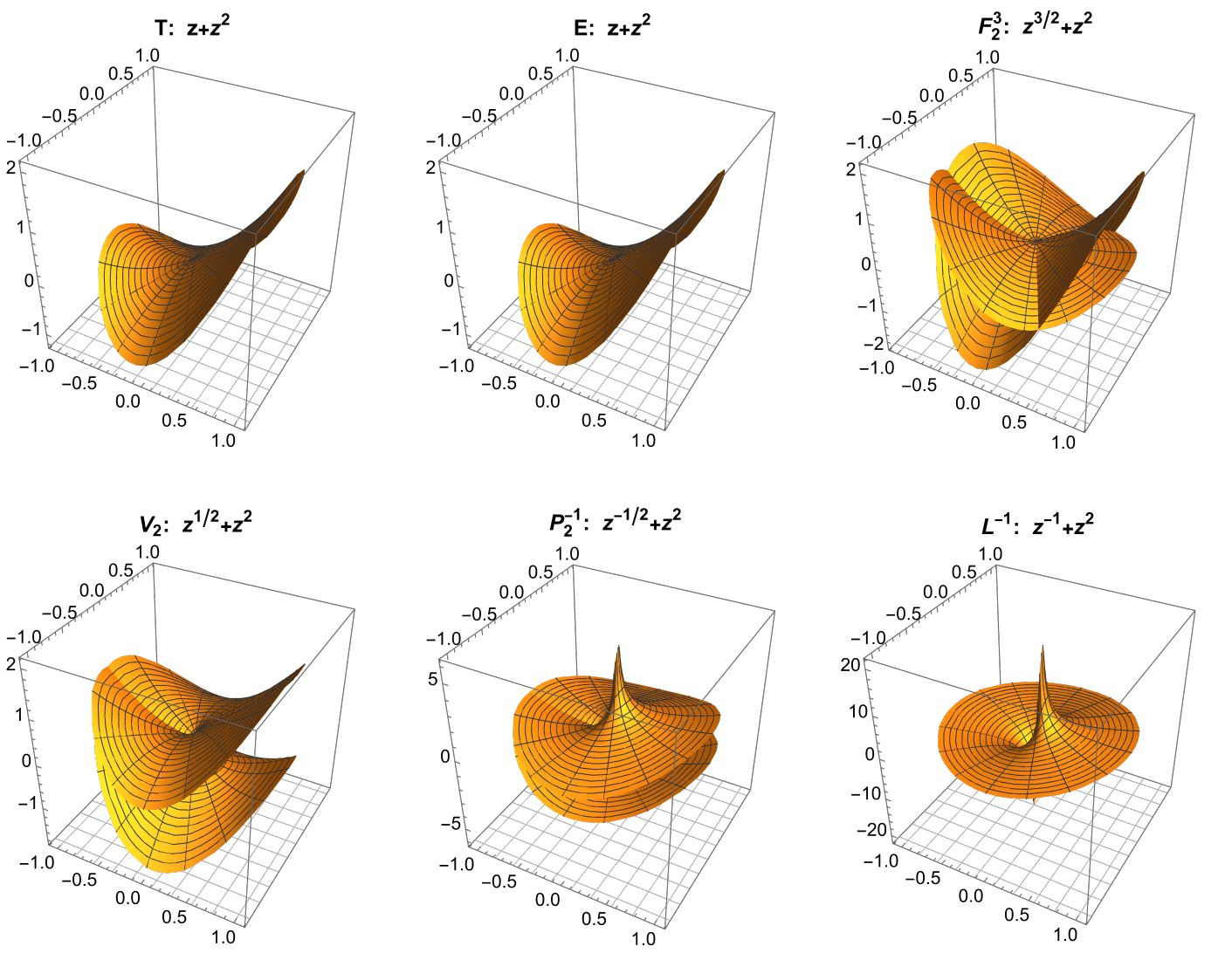}
\caption{Branch Types (real surfaces)}
\label{figure:appendixa}
\end{figure}
\pagebreak

\section{}
\begin{center}
Details about Puiseux expansions of algebraic functions
\label{appendix:appB}
\end{center}
Given (\ref{eqn:eqn001}), $w(z)$ is expanded at $z_0$ into a set of $n$ series segregated into sets of conjugate cycles $C=\{c_i\}$ given by
\begin{equation}
 \{P_i(z)\}=\{w_i(z)\}=\left\{\sum_{k=r}^\infty a_k \left((z-z_0)^{\frac{1}{c_i}}\right)^{m_k}\right\};\quad, i=1,2\cdots,C.
\label{app2eqn002}
\end{equation}
with $r$ related to the series order, and the set $\{m_k\}$, exponent numerators over a least common denominator $c_i$.
Each set of cycles $c_i$ contains $c_i$ series since 
\begin{equation}
(z-z_0)^{\frac{m_k}{c_i}}\equiv |r|^{\frac{m_k}{c_i}}\left(e^{\frac{i}{c_i}(\theta+2u\pi)}\right)^{m_k},\quad u=0,1,\cdots,c_i-1.
\end{equation}

 For example, if a $10$ degree function ramified into a $1,2,3$ and $4$ cycle branch, the set of cycle indexes would be $\{c_i\}=\{1,2,3,4\}$.   In this case, $\{P_i(z)\}$ is a set of ten series consisting of four cycles.  The first cycle set contains one $1$-cycle branch, the second set contains two $2$-cycle series, the third, three $3$-cycle series, and the last, four $4$-cycle series.    However, Puiseux series are easier to work with if $(z-z_0)\to z$ in (\ref{app2eqn002}):
\begin{equation}
 \{P_j(z)\}=\{w_j(z)\}=\left\{\sum_{k=r}^\infty a_k z^{\frac{m_k}{c_j}}\right\}_{j=1}^C;\quad z\to (z-z_0)
\label{app2eqn002b}
\end{equation} 

with a particular value of the expansion at $t$ given by $P_j(t-z_0)$. 
 
  And when Puiseux series are generated by Newton Polygon as is done in this paper, each of the $c_i$ series in each conjugate class of series are already in conjugate form;  the $\displaystyle \left(e^{\frac{2u\pi i}{c_j}}\right)^{m_k}$ terms have been absorbed into the coefficients producing $n$ series:
	\begin{equation}
 \{P_i(z)\}=\{w_i(z)\}=\left\{\sum_{k=r}^\infty a_k z^{\frac{m_k}{c_i}}\right\}_{i=1}^n;z^{\frac{m_k}{c_i}}\equiv |r|^{\frac{m_k}{c_i}}\left(e^{\frac{i \theta}{c_i}}\right)^{m_k}
\label{appeqn002c}
\end{equation} 

with $\displaystyle z^{\frac{m_k}{c_i}}$ interpreted as principal-valued and $c_i$ the cycle size.  And Puiseux series are easier to understand initially by simply viewing Puiseux series one at a time:
\begin{equation}
 P(z)=\sum_{k=r}^\infty a_k z^{\frac{m_k}{c}}
\label{eqn:eqn002d}
\end{equation} 
where $z^{\frac{m_k}{c}}$ is interpreted as principal-value, a particular value $t$ of the series is given by $P(t-z_0)$ where $z_0$ is the expansion center, the set $\{m_i\}$ numerators of exponents placed under a lowest common denominator $c$,  and the full set of Puiseux expansions at $z_0$ is $n$ such series that can have different values of $c$.
\pagebreak

\end{document}